\documentclass[letterpaper]{JHEP3}
\usepackage{epsfig,amsmath,xspace}

\newcommand{\Z}{\mathbb{Z}}
\newcommand{\R}{\mathbb{R}}

\newcommand{\pd}[2]{\frac{\partial #1}{\partial #2}}
\newcommand{\K}{K\"ahler\xspace}
\newcommand{\CP}{\mathbb{CP}}

\title{Numerical \K-Ricci soliton on the second del Pezzo}

\author{Matthew Headrick \\ Stanford Institute for Theoretical Physics, Stanford, CA 94305-4060, USA \\ \email{headrick@stanford.edu}}

\author{Toby Wiseman \\ Blackett Laboratory, Imperial College London, London SW7 2AZ, UK \\ \email{t.wiseman@ic.ac.uk} }

\abstract{The second del Pezzo surface is known by work of Tian-Zhu and Wang-Zhu to admit a unique \K-Ricci soliton. Applying a method described in hep-th/0703057, we use Ricci flow to numerically compute that soliton metric. We numerically compute the value of its Perelman entropy (or Gaussian density).}

\preprint{SU-ITP-07/07 \\ Imperial/TP/07/TW/02}

\begin{document}

In a recent paper \cite{Doran}, Doran, Herzog, Kantor, and the present authors presented several numerical methods for solving the Einstein equation on toric manifolds. We applied those methods to the third del Pezzo surface ($\CP^2\#3\overline{\CP^2}$, or $dP_3$), which is known by work of Tian-Yau \cite{MR894378, MR904143} and Siu \cite{MR942521} to admit a unique \K-Einstein metric. Two of the numerical methods involved simulation of (normalized) Ricci flow,
\begin{equation}\label{RFdef}
\pd{g_{\mu\nu}}{t} = -2R_{\mu\nu} + 2g_{\mu\nu},
\end{equation}
which is guaranteed by a theorem of Tian-Zhu \cite{TianZhuConvergence} to converge starting from any \K metric in the same class as the Ricci form.

In this note we show that simulation of Ricci flow can also be applied effectively to  toric manifolds that do not admit a \K-Einstein metric but rather a \emph{\K-Ricci soliton}.\footnote{See \cite{MR2243675} for a review of Ricci solitons.} A (shrinking) \K-Ricci soliton consists of a \K metric $g_{\mu\nu}$ and a holomorphic vector field $\xi^\mu$ satisfying
\begin{equation}\label{solitondef}
R_{\mu\nu} = g_{\mu\nu} - \nabla_{(\mu}\xi_{\nu)}.
\end{equation}
Such a metric is a fixed point up to diffeomorphisms of the flow \eqref{RFdef}, and the Tian-Zhu theorem again guarantees convergence (up to diffeomorphisms) to it. In this paper we focus specifically on the second del Pezzo surface ($dP_2$). The existence and uniqueness of a \K-Ricci soliton on $dP_2$ follow from theorems of Wang-Zhu \cite{MR2084775} and Tian-Zhu \cite{MR1768112, MR1915043} respectively. However, this metric is not known explicitly. We will describe a numerical approximation to it obtained using one of the methods of \cite{Doran}, and explore some of its geometrical properties. (The first del Pezzo surface is also toric, and, like $dP_2$, admits a \K-Ricci soliton \cite{MR1145263}; however, that soliton is co-homogeneity 1 and is already known in a fairly explicit form \cite{MR1417944}.)

We begin by reviewing some basic facts about toric metrics and $dP_2$. A toric manifold admits two natural coordinate systems: complex coordinates $u^i+i\theta^i$, and symplectic coordinates $(x_i,\theta^i)$. The $\theta^i$ are periodic ($\theta^i\sim\theta^i+2\pi$), and the $U(1)^n$ isometries (where $n$ is the complex dimension) of any toric metric act by translating them, leaving $u^i$ and $x_i$ fixed. The real parts $u^i$ of the complex coordinates range over $\R^n$, and the metric can be expressed in terms of the \K potential $f(u)$ as
\begin{equation}
ds^2 = F_{ij}(du^idu^j+d\theta^i d\theta^j),
\end{equation}
where
\begin{equation}
F_{ij} = \frac{\partial^2 f}{\partial u^i\partial u^j}.
\end{equation}
The symplectic coordinates $x_i$ range inside a Delzant polytope, which in the case of $dP_2$ is a pentagon, and the metric can be expressed in terms of the symplectic potential $g(x)$ as
\begin{equation}\label{sympmetric}
ds^2 = G^{ij}dx_idx_j + G_{ij}d\theta^i d\theta^j,
\end{equation}
where
\begin{equation}
G^{ij} = \frac{\partial^2 g}{\partial x_i\partial x_j}
\end{equation}
and $G_{ij}$ is the inverse matrix of $G^{ij}$. The \K and symplectic potentials are related by a Legendre transform,
\begin{equation}
x_i = \pd{f}{u^i}, \qquad u^i = \pd{g}{x_i}, \qquad f + g = u^ix_i,
\end{equation}
so the mapping between $u$ and $x$ depends on the metric. Under this mapping $F_{ij} = G_{ij}$. The \K class is related to the positions of the edges of the Delzant polytope. For $dP_2$, we are interested in the case where the \K class is equal to the class of the Ricci tensor ($[\omega]=2\pi c_1$); the corresponding polytope is as shown in figure \ref{polytope}. The Ricci tensor is given by
\begin{equation}
R_{u^iu^j} = R_{\theta^i\theta^j} = \frac{\partial^2r}{\partial u^i\partial u^j}
\end{equation}
where
\begin{equation}
r = -\frac12\log\det(F_{ij}).
\end{equation}

\FIGURE{
\epsfig{file=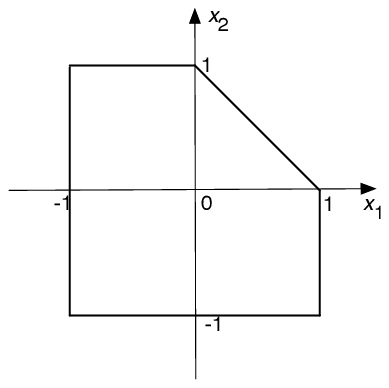,width=2.7in}
\caption{The Delzant polytope for $dP_2$, with $[\omega]=2\pi c_1$.}
\label{polytope}
}

Several facts about the \K-Ricci soliton on $dP_2$ can be deduced analytically. First, the Wang-Zhu existence theorem implies that it is a \emph{gradient} soliton, meaning that $\xi_\mu=\partial_\mu\phi$ for some (globally defined) function $\phi$. (This also follows from a theorem by Perelman that any compact shrinking Ricci soliton is a gradient soliton \cite{perelman-2002}.) It follows that ${J^\mu}_\nu\xi^\nu$, where ${J^\mu}_\nu$ is the manifold's complex structure, is a holomorphic Killing field. (From the symplectic viewpoint, $\phi$ is a Hamiltonian function for the holomorphic isometry generated by ${J^\mu}_\nu\xi^\nu$.) Assuming that the metric has no continuous isometries other than its toric ones, $\xi^\mu\partial_\mu$ must be a linear combination of $\partial_{u^1}$ and $\partial_{u^2}$. Assuming further that the polytope's $\Z_2$ symmetry under which $x_1$ and $x_2$ are exchanged (which lifts to the diffeomorphism $u^1+i\theta^1 \leftrightarrow u^2+i\theta^2$ of the manifold) is preserved by the soliton, we learn that
\begin{equation}\label{phiform}
\xi_{x_i} = \xi^{u^i} = \alpha, \qquad \xi^{\theta^i} = 0,
\qquad \phi = \alpha(x_1+x_2)
\end{equation}
where $\alpha$ is a constant. We will see that the numerics bear out these assumptions.

We simulated Ricci flow on $dP_2$ using a method analogous to that described in subsection 4.2 of \cite{Doran}. (Details of the implementation were analogous to those described in appendix B of that paper, with a maximum lattice resolution of $512\times512$. The code can be downloaded from the website \cite{website}.) The initial metric was taken to be the ``canonical" one of Guillemin \cite{MR1293656} and Abreu \cite{MR1969265}, with symplectic potential
\begin{equation}
g_{\rm can} = \frac12\sum_{a=-2}^2l_a\log l_a,
\end{equation}
where
\begin{equation}
l_0 = 1-(x_1+x_2), \qquad l_{\pm1} = 1\pm x_1, \qquad l_{\pm2} = 1\pm x_2.
\end{equation}
At late times, the metric converged to one whose symplectic potential $g$ satisfied the equation
\begin{equation}
r+g-x_i\pd{g}{x_i} = -\frac\alpha2(x_1+x_2),
\end{equation}
implying equation \eqref{solitondef}, where
\begin{equation}
\alpha \approx -0.43475.
\end{equation}
The numerical data encoding this symplectic potential, along with a Mathematica notebook for manipulating that data, are available on the website \cite{website}. However, a good approximation may be obtained by fitting the difference between $g$ and $g_{\rm can}$ to quartic order in $x_i$, yielding
\begin{equation}\begin{split}
g \approx&\,g_{\rm can} \\
& - 0.087 x_1 x_2 - 0.121 (x_1^2+x_2^2) \\
& - 0.041 x_1x_2(x_1+x_2) - 0.031(x_1^3+x_2^3) \\
& - 0.015 x_1^2x_2^2 - 0.013 x_1x_2(x_1^2+x_2^2) -0.009(x_1^4+x_2^4)
.
\end{split}\end{equation}
This yields metric components $G^{ij}$ with an absolute error of less than $0.02$. Three curvature invariants for the soliton metric are plotted in figure \ref{curvatures}.

\FIGURE{
\epsfig{file=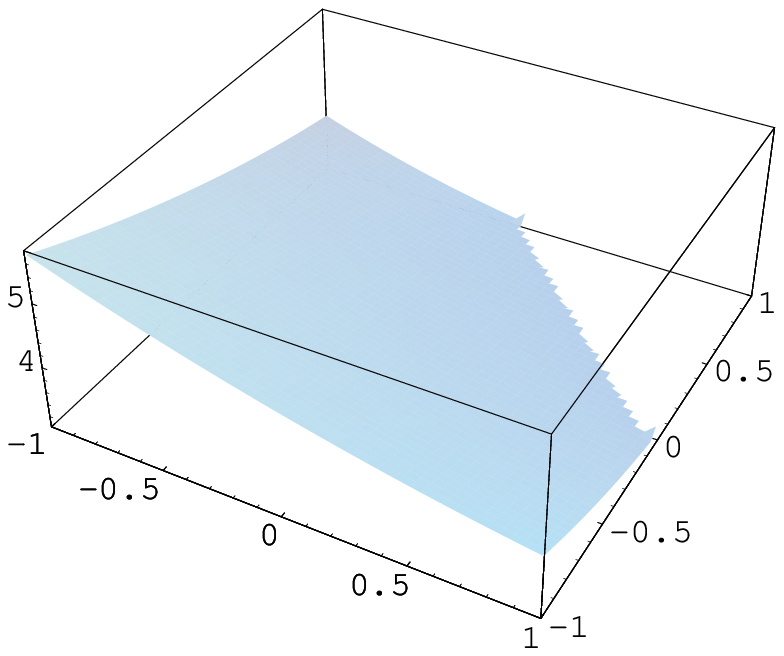,width=2.7in}\\
\epsfig{file=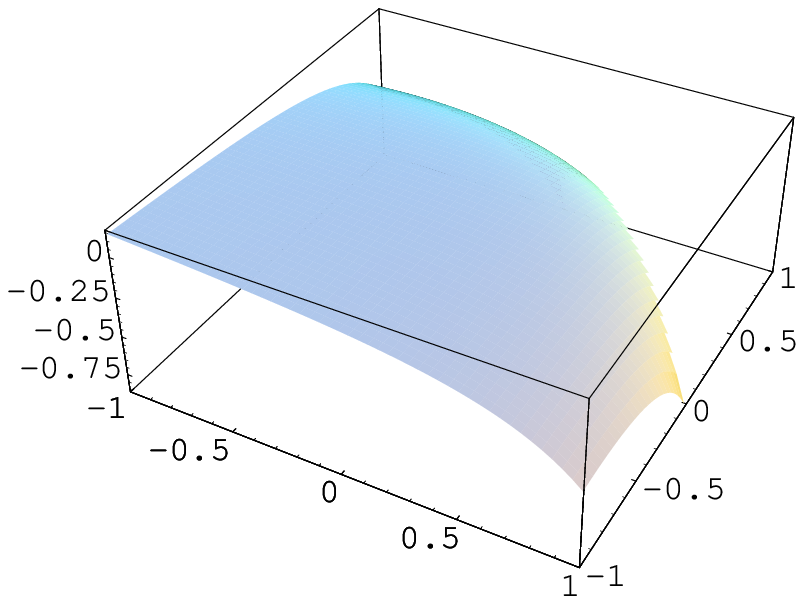,width=2.7in}\\
\epsfig{file=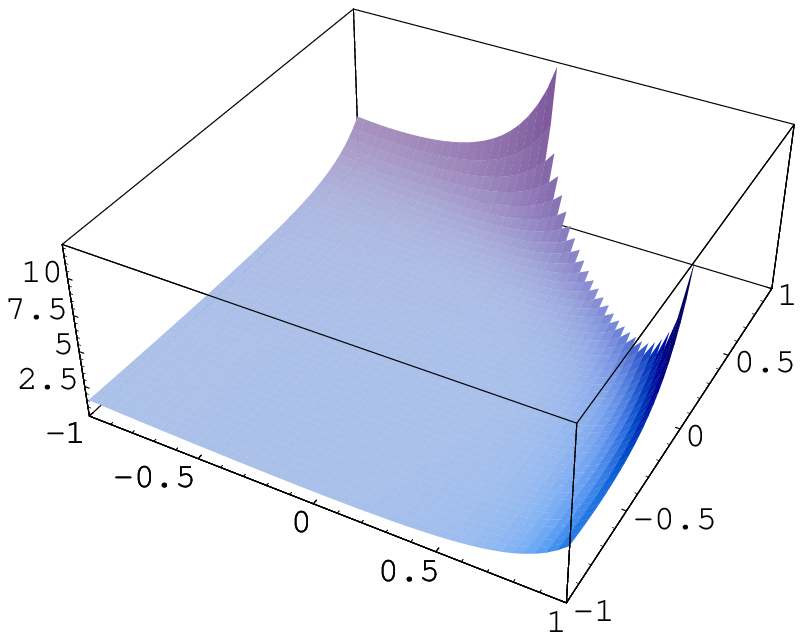,width=2.7in}\\
\caption{
Three curvature invariants for the Ricci soliton metric, versus $x_1$ and $x_2$: (top) Ricci scalar; (middle) $R_{x_1x_2x_1x_2}/(G^{11}G^{22}-(G^{12})^2)$, the sectional curvature of a surface of constant $\theta$; (bottom) $(R_{\mu\nu\kappa\lambda}R^{\mu\nu\kappa\lambda}-4R_{\mu\nu}R^{\mu\nu}+R^2)/8$ ($4\pi^2$ times the Euler density), which integrates to 5, the Euler character of $dP_2$.
}
\label{curvatures}
}

Cao \cite{MR2243675} lists the calculation of Perelman's entropy functional \cite{perelman-2002} for the soliton on $dP_2$ as an important open problem. This is straightforward to compute given our numerical solution. For any gradient Ricci soliton, the Perelman entropy is the standard entropy (at unit temperature) of the canonical ensemble defined by taking the manifold as a phase space and $\phi$ as the Hamiltonian. Specifically, we define the partition function and entropy as usual by
\begin{equation}
Z(\beta) = \frac1{(2\pi e)^n}\int\sqrt{g}e^{-\beta\phi}, \qquad
S(\beta) = \left(1-\beta\frac{d}{d\beta}\right)\log Z(\beta)
\end{equation}
(the pre-factor of $(2\pi e)^{-n}$ in $Z$ is purely conventional). The Perelman entropy is then $\nu = S(1)$. Specializing to the toric case, since in symplectic coordinates $\sqrt{g}=1$ (see \eqref{sympmetric}), we see that to compute $\nu$ it is only necessary to know $\phi$, not the metric itself. On $dP_2$, for $\phi$ as in equation \eqref{phiform}, we find
\begin{eqnarray}
Z(\beta) &=& 
\frac1{(\beta\alpha e)^2}\left((1-\beta\alpha)e^{-\beta\alpha}-2+e^{2\beta\alpha}\right), \\
\nu &=& \frac{\alpha  \left(1+2 e^{\alpha }-3 e^{3 \alpha
   }\right)}{1-\alpha -2 e^{\alpha }+e^{3 \alpha }}+\log
   \left(\frac{1-\alpha -2 e^{\alpha }+e^{3 \alpha
   }}{\alpha ^2}\right)
\approx -0.78778
\end{eqnarray}
Cao-Hamilton-Ilmanen \cite{cao-2004} define the \emph{Gaussian density} of a soliton as $\Theta=e^\nu$; we have $\Theta \approx 0.45485 \approx 3.3609/e^2$. This value fits nicely into the pattern of values observed for other \K-Ricci solitons on del Pezzo surfaces, lying between those for the Koiso soliton on $dP_1$ ($3.826/e^2$) and for the \K-Einstein metric on $dP_3$ ($3/e^2$).

\section*{Acknowledgments}

We would like to thank C. Doran, C. Herzog, and J. Kantor for collaboration and discussions in the early stages of this project, and especially C. Doran for suggesting this problem. We would also like to thank D. Knopf and D. Martelli for useful conversations.
M.H. would like to thank the Theory Group of the University of Texas, Austin and the MIT Center for Theoretical Physics for hospitality while this work was being done. He is supported  by the Stanford Institute for Theoretical Physics and by NSF grant PHY 9870115.
T.W. is supported by a PPARC advanced fellowship and the Halliday award.

%
\bibliography{ref}

\providecommand{\href}[2]{#2}\begingroup\raggedright\begin{thebibliography}{10}

\bibitem{Doran}
C.~Doran, M.~Headrick, C.~P. Herzog, J.~Kantor, and T.~Wiseman, {\it Numerical
  {K}aehler-{E}instein metric on the third del {P}ezzo},
  \href{http://xxx.lanl.gov/abs/hep-th/0703057}{{\tt hep-th/0703057}}.

\bibitem{MR894378}
G.~Tian, {\it On {K}\"ahler-{E}instein metrics on certain {K}\"ahler manifolds
  with {$C\sb 1(M)>0$}},  {\em Invent. Math.} {\bf 89} (1987), no.~2 225--246.

\bibitem{MR904143}
G.~Tian and S.-T. Yau, {\it K\"ahler-{E}instein metrics on complex surfaces
  with {$C\sb 1>0$}},  {\em Comm. Math. Phys.} {\bf 112} (1987), no.~1
  175--203.

\bibitem{MR942521}
Y.~T. Siu, {\it The existence of {K}\"ahler-{E}instein metrics on manifolds
  with positive anticanonical line bundle and a suitable finite symmetry
  group},  {\em Ann. of Math. (2)} {\bf 127} (1988), no.~3 585--627.

\bibitem{TianZhuConvergence}
G.~Tian and X.~Zhu, {\it Convergence of {K}\"ahler-{R}icci flow},  {\em J.
  Amer. Math. Soc.} {\bf 20} (2007) 675--699.

\bibitem{MR2243675}
H.-D. Cao, {\it Geometry of {R}icci solitons},  {\em Chinese Ann. Math. Ser. B}
  {\bf 27} (2006), no.~2 121--142.

\bibitem{MR2084775}
X.-J. Wang and X.~Zhu, {\it K\"ahler-{R}icci solitons on toric manifolds with
  positive first {C}hern class},  {\em Adv. Math.} {\bf 188} (2004), no.~1
  87--103.

\bibitem{MR1768112}
G.~Tian and X.~Zhu, {\it Uniqueness of {K}\"ahler-{R}icci solitons},  {\em Acta
  Math.} {\bf 184} (2000), no.~2 271--305.

\bibitem{MR1915043}
G.~Tian and X.~Zhu, {\it A new holomorphic invariant and uniqueness of
  {K}\"ahler-{R}icci solitons},  {\em Comment. Math. Helv.} {\bf 77} (2002),
  no.~2 297--325.

\bibitem{MR1145263}
N.~Koiso, {\it On rotationally symmetric {H}amilton's equation for
  {K}\"ahler-{E}instein metrics},  in {\em Recent topics in differential and
  analytic geometry}, vol.~18 of {\em Adv. Stud. Pure Math.}, pp.~327--337.
\newblock Academic Press, Boston, MA, 1990.

\bibitem{MR1417944}
H.-D. Cao, {\it Existence of gradient {K}\"ahler-{R}icci solitons},  in {\em
  Elliptic and parabolic methods in geometry (Minneapolis, MN, 1994)},
  pp.~1--16.
\newblock A K Peters, Wellesley, MA, 1996.

\bibitem{perelman-2002}
G.~Perelman, {\it The entropy formula for the {R}icci flow and its geometric
  applications},  \href{http://xxx.lanl.gov/abs/math/0211159}{{\tt
  math/0211159}}.

\bibitem{website}
http://www.stanford.edu/\~{}headrick/dP2/.

\bibitem{MR1293656}
V.~Guillemin, {\it Kaehler structures on toric varieties},  {\em J.
  Differential Geom.} {\bf 40} (1994), no.~2 285--309.

\bibitem{MR1969265}
M.~Abreu, {\it K\"ahler geometry of toric manifolds in symplectic coordinates},
   in {\em Symplectic and contact topology: interactions and perspectives
  (Toronto, ON/Montreal, QC, 2001)}, vol.~35 of {\em Fields Inst. Commun.},
  pp.~1--24.
\newblock Amer. Math. Soc., Providence, RI, 2003.

\bibitem{cao-2004}
H.-D. Cao, R.~S. Hamilton, and T.~Ilmanen, {\it Gaussian densities and
  stability for some {R}icci solitons},
  \href{http://xxx.lanl.gov/abs/math/0404165}{{\tt math/0404165}}.

\end{thebibliography}\endgroup
\bibliographystyle{JHEP}
%

\end{document}